\documentclass[11pt]{article}

\usepackage[margin=1in]{geometry}
\usepackage{amsmath,amssymb,amsfonts,mathtools}
\usepackage{siunitx}
\usepackage{booktabs,multirow}
\usepackage{graphicx}
\usepackage{soul}
\usepackage{xcolor}
\usepackage{subcaption}
\usepackage{lineno}
\usepackage{booktabs}
\usepackage{siunitx}
\usepackage{acronym}
\usepackage{multirow}
\usepackage{hyperref}
\usepackage[capitalize,noabbrev]{cleveref}
\usepackage{algorithm}
\usepackage{authblk}
\usepackage{setspace}
\usepackage{verbatim}
\usepackage{algpseudocode}
\usepackage[backend=biber,style=ieee,maxbibnames=6]{biblatex}
\usepackage{xcolor}
\newcommand{\DarkComment}[1]{\textcolor{darkgreen}{#1}}

\definecolor{darkgreen}{RGB}{0,100,0}

\algrenewcommand{\algorithmiccomment}[1]{\hfill\DarkComment{\(\triangleright\) #1}}

\allowdisplaybreaks[4]

\acrodef{MIP}[MIP]{Mixed-Integer Program}
\acrodef{UAS}[UAS]{Uncrewed Aerial Systems}
\acrodef{CoVRP}[CoVRP]{Coverage Vehicle Routing Problem}
\acrodef{VRP}[VRP]{Vehicle Routing Problem}
\acrodef{TSP}[TSP]{Traveling Salesman Problem}
\acrodef{MCLP}[MCLP]{Maximal Covering Location Problem}

\title{\textbf{Optimizing Agricultural Drone Operations:\\From Launch and Recovery Siting to\\Tiered Routing Strategies}}

\author[1]{Ethan Kolby\footnote{PhD Student, Department of Aerospace Engineering, University of Michigan, Ann Arbor, Michigan 48109.}}
\author[2]{Josh Noble\footnote{Principal Data Scientist, Mosaic ATM.}}
\author[1]{Max Z. Li\footnote{Assistant Professor, Department of Aerospace Engineering, Department of Civil and Environmental Engineering, Department of Industrial and Operations Engineering, University of Michigan, Ann Arbor, Michigan 48109, AIAA Member.}}
\affil[1]{University of Michigan}
\affil[2]{Mosaic ATM}

\begin{document}

\maketitle

\begin{abstract}

Drones are increasingly used in agriculture, where tight margins demand efficient planning. Current optimization tools suffer from exponential runtimes as problem sizes grow, necessitating practical heuristics for daily operations. This paper presents an operational framework and benchmarking analysis for drone spraying operations. We evaluate the trade-offs between facility siting methods and tiered routing parameters. For facility siting, comparing a Mixed-Integer Program (MIP) baseline against a $p$-Median heuristic shows that the heuristic reduces runtime by three orders of magnitude, from over 97 seconds to under 1.2 seconds, with only a 4\% reduction in serviced field area. For route planning, a tiered problem decomposition approach partitioning the target area into 6 to 8 spatial clusters reduces computation time by an order of magnitude with minimal degradation in serviced area. This framework achieves minute-scale planning on commodity hardware, demonstrating operational relevance. Future research will incorporate weather modeling, integrated optimization of facility location and routing, and validation across diverse field geometries.

\end{abstract}
\section{Introduction}
\label{sec:intro}

Precision agriculture is increasingly relying on \ac{UAS} to optimize crop management and chemical application. For operations like spot-spraying, drones offer unparalleled agility compared to traditional heavy machinery. However, this sector operates on exceptionally tight margins. Drone deployments are bounded by typical endurance limits of 15 to 30 minutes and limited payload weight. To make these systems economically viable, operators require highly efficient, minute-scale mission planning that maximizes the area sprayed per battery charge while minimizing non-productive transit time.

This operational challenge naturally frames crop spraying as a \ac{CoVRP}. Unlike the classical \ac{VRP} that serves discrete point locations, the \ac{CoVRP} requires a vehicle to satisfy continuous regional targets via a specific footprint, such as a spray radius. Feasibility and optimality depend on how flight paths sweep through regions, managing boundary constraints and overlapping passes. Unfortunately, existing optimal decision-making tools and exact mathematical formulations for the \ac{CoVRP} suffer from exponentially increasing runtimes as field sizes and geometric complexity grow. This computational bottleneck renders exact solvers impractical for a farmer or operator standing at the edge of a field waiting to launch a fleet.

To bridge this gap between complex coverage requirements and strict operational time constraints, this paper introduces an applied operational framework that relies on hierarchical decomposition. The pipeline begins with facility selection, where we adapt a greedy $p$-Median heuristic to act as a rapid spatial partitioner. To evaluate this approach, we compare the heuristic against a load-balancing optimization method to explicitly examine the trade-off between initial computational speed and downstream routing performance. Once facilities are sited, the framework further decomposes the assigned areas by grouping the routing nodes using $k$-means clustering. A high-level \ac{TSP} is then solved to establish the optimal sequence of these subgroups. Finally, a discrete \ac{CoVRP} is solved within each individual subgroup before the segments are stitched together into a complete, battery-constrained flight route. By sequentially breaking the problem down, this end-to-end architecture generates near-optimal field coverage in minutes rather than hours.

While precision agriculture serves as the primary testbed for this benchmarking analysis, the structural elements of this framework are highly adaptable. Sectors with tight time and energy budgets face similar operational hurdles: last-mile delivery, infrastructure inspection, and emergency response all share aspects of the core \ac{CoVRP} architecture. Whether a drone is dropping parcels on campuses, using photogrammetry to create 3D models of bridges \cite{congress2024trr}, or deploying sensors for real-time situational awareness during wildfires and search-and-rescue missions \cite{ragbir2023sensors,lyu2023remotesensing}, the underlying need to rapidly pair facility siting with capacity-constrained spatial coverage remains a universal challenge.

\subsection{Background and Motivation}
\label{sec:background}

The classical \ac{VRP} serves discrete customers at points, with feasibility defined by visiting each node once under vehicle and depot constraints. The \ac{CoVRP} replaces point service with regional service: a target is satisfied when the vehicle’s footprint (extent of spray coverage, sensor field of view, or communication radius) intersects its associated region \cite{glock2023spatial}. This change couples geometry and routing. Feasibility depends on how paths pass through regions, not just on node visitation, and performance must account for coverage quality near boundaries and in overlaps. As a result, optimality involves more than distance minimization; it also involves covered and uncovered area, as well as redundant passes.

Operations across delivery, inspection, emergency response, and agriculture share the same structural elements. The area of interest is discretized into coverage nodes, where each node has a coverage neighborhood that is served when the node is visited. Edges between nodes represent travel costs, quantified here as distance. A subset of candidate facilities is activated as \ac{UAS} launch and recovery sites, using a facility location planning algorithm. Vehicles have finite range, payload, and kinematic limits, and the resulting routes must satisfy both coverage and distance constraints. Poor siting can lead to coverage gaps if there are no selected facilities within the maximum distance. Furthermore, siting determines the workload distribution across vehicles, potentially increasing the overall makespan if one facility is overloaded. These considerations motivate a direct comparison between a load-balancing siting optimizer and a greedy location siting heuristic for a fixed facility budget.

We have identified siting strategies, tiered routing architectures, and adaptable coverage models as the primary operational levers for this problem. To set up the subsequent head-to-head experiments, Section~\ref{sec:lit_review} situates our work within the literature on spatial coverage, ``close-enough'' \ac{VRP} formulations \cite{mennell2009heuristics}, \ac{UAS} facility location, and hierarchical solution pipelines.

\subsection{Literature Review}
\label{sec:lit_review}

\subsubsection{Coverage and Location-Routing}

Spatial-coverage formulations generalize point-service \ac{VRP}s by satisfying targets when a vehicle's footprint intersects a tolerance region. A structured classification of spatial coverage models is crucial, as it clarifies core architectural decisions and reveals established solution trends, such as decomposition and large-neighborhood search. \cite{glock2023spatial}. However, the literature has historically lacked a systematic investigation of these decomposition strategies, a gap addressed by Santini et al., who find that route-based decompositions are generally superior to path-based methods \cite{santini2023decomposition}.

Two foundational streams underlie modern \ac{CoVRP} work. First, the idea of surrogate visitation originates with the Covering Salesman Problem (CSP) \cite{current1989csp} and the Covering Tour Problem (CTP) \cite{gendreau1997ctp}, which combine \ac{TSP} style routing with coverage constraints. These seminal papers formalize Integer Linear Programming (ILP) models and introduce heuristics that remain the basis for many modern relaxations \cite{gendreau1997ctp}. Second, the ``close-enough" family of problems (e.g., CETSP) sharpens the geometry of these tolerance regions and provides exact and approximate solvers \cite{behdani2014cetsp, coutinho2016cetsp}. These foundational results establish that covering constraints are natural surrogates for footprint-based tasks like sensing or spraying and that scalable solutions typically rely on hierarchical decompositions. For example, Cao et al. demonstrate a hierarchical planner for complex 3D environments that first subdivides the space and then solves local TSPs, resulting in a tenfold speedup over non-decomposed methods \cite{cao2020hierarchical}.

When facility location siting and vehicle routing are co-decided, the location–routing viewpoint becomes essential. The classic CTP already illustrates how the selection of visited nodes trades off against tour length, motivating two-stage pipelines that first pick sites or depots and then routes vehicles \cite{gendreau1997ctp}. This trade-off is central to the Maximum Coverage Capacitated Facility Location Problem with Drones (MCFLPD), where Chauhan et al. explicitly quantify the downstream coverage loss at nearly 20\% when substituting a fast greedy siting heuristic for an optimal MIP-based approach \cite{chauhan2019mcflpd}. This hierarchical structure appears vividly in hybrid ground–air systems. Murray and Chu’s Flying-Sidekick TSP, for example, formalizes synchronized drone launches from a truck, proving the value of co-optimizing the truck route and aerial sorties \cite{murray2015fstsp}. More recent work reinforces this structure by co-optimizing depot selection, synchronization, and last-mile tours in collaborative truck-drone delivery, yielding substantial mission-time savings \cite{springer_robot_stations_2023, zieher2024parcel, fu2025truckdrone, zhen2025truckdrone}. Within these hierarchies, classic local-search operators like 2-opt and Lin–Kernighan are mainstays for fast route refinement \cite{croes1958twoopt, lin1973lk}.

\subsubsection{Operational Constraints in Agricultural Spraying}

The specific constraints and objectives of real-world applications guide the design of CoVRP algorithms and metrics. Unlike the standard \ac{VRP} which minimizes distance to visit a set of discrete nodes, CoVRP variants must maximize a continuous or discretized coverage metric—such as sprayed area, gathered sensor data, or probability of detection—under strict resource constraints. Endurance limitations are a common thread, as contemporary multirotors typically sustain only 20–30 minutes of flight, compared to hour-class flights for fixed-wing platforms. This necessitates a shift from simple ``lawnmower" patterns to sophisticated, tiered routing choices where depot placement and area selection are coupled \cite{peksa2024sensors, silberstein2025jaot, mitridis2023aerospace}. 

Recent surveys document the rapid uptake of this tiered approach for agricultural \ac{UAS} tasks like spot-spraying, where optimization frameworks must balance the granular ``value" of covering a specific node against the battery cost of reaching it \cite{tsouros2019review, mdpi_weed_mip_2024}. For example, deep learning models now generate ``weed pressure maps," converting the field into a weighted grid; the CoVRP algorithm must then determine which high-weight cells can be visited within a strict energy budget \cite{deep_learning_weed_2023}. Crucially, agricultural environments introduce geometric complexities that standard routing ignores. Wang et al. apply a logic-based optimization algorithm to determine optimal spray headings in irregular fields, aiming to minimize the ``non-productive" flight time spent turning at the field boundaries \cite{wang2025ddpg}. Their work highlights that in CoVRP, the quality of coverage is as important as the quantity. However, existing agricultural literature often assumes a static depot location at the field edge. This motivates our investigation into optimal facility placement; by optimizing the launch point relative to the high-value infestation clusters, we can significantly increase the reachable coverage area per battery cycle.

\subsubsection{Sensor-Based Coverage}

While precision agriculture focuses on the physical application of agrichemicals, the broader CoVRP framework is highly adaptable to other domains where ``coverage'' is achieved passively via onboard sensors and cameras. In these applications, the structural imperative remains the same: maximizing a continuous coverage metric under vehicle constraints. 

For instance, in infrastructure inspection, \ac{UAS} equipped with photogrammetry or LiDAR must satisfy specific visual constraints, such as incidence angles, pixel resolution, and image overlap, to create high-fidelity 3D digital twins \cite{congress2024trr, gaspari2022isprs, zollini2020remotesensing}. Similarly, in Search-and-Rescue (SAR) operations, the objective shifts to maximizing the ``Probability of Detection''. This transforms the routing problem into a submodular maximization task where the marginal gain of sensor coverage diminishes as an area is repeatedly scanned \cite{waharte2010sar, lyu2023remotesensing, zhang2025jrs, shaheen2023access}. Finally, wildland fire management utilizes sensor-based coverage for real-time monitoring under extreme temporal constraints. These missions demand rapid, heuristic-driven routing that can account for intense convective winds and shifting burn windows \cite{ragbir2023sensors, bouguettaya2022sigpro, lopresti2024ijwf}. 

Although the specific definition of coverage (e.g., visual viewpoints, detection probability, or fire monitoring) differs from crop spraying, these domains all share the fundamental need for scalable, tiered routing architectures that closely couple facility placement with the routing of irregular coverage zones.
\subsection{Technical Gaps}
\label{sec:gaps}

Despite strong baselines for routing problems and growing interest in spatial coverage, several gaps constrain the deployment of practical coverage routing solutions for applications like crop spraying.

First, while head-to-head comparisons of facility siting methods (e.g., optimal \ac{MIP}s versus greedy heuristics) do exist, the literature has not sufficiently explored how these upstream choices propagate through multi-stage planning pipelines to affect downstream performance. Prior studies often evaluate siting strategies in isolation. They do not quantify how substituting an optimal but slow facility selection method for a fast heuristic one alters final tour metrics like total coverage, makespan, and load balance after subsequent routing and optimization stages are completed. This leaves a gap in understanding the full, end-to-end trade-offs between optimality and speed in a tiered system.

Second, while hierarchical or ``cluster-first, route-second" approaches are known to be effective for large-scale VRPs, the literature seldom maps their joint performance–runtime frontiers for coverage objectives. These tiered pipelines expose highly interdependent parameters—such as the number of clusters, the method for sequencing clusters, and the time limits allocated to subproblems—yet most reports fix these knobs and emphasize algorithmic novelty. As recent surveys have noted, the literature lacks a systematic investigation of how different decomposition strategies impact state-of-the-art heuristics \cite{santini2023decomposition}. This leaves a practical gap in guidance for tuning a tiered CoVRP system to a specific fleet, field geometry, and compute budget.

Finally, while wall-clock performance is often reported, achieving minute-scale or faster runtimes on large, complex, and diverse geometries remains a critical challenge. Many studies emphasize solution quality on modest or standardized instance sizes. However, for operational deployment where plans must be recomputed quickly as conditions or assets change, the ability to generate high-quality solutions within a tight time budget is paramount. The performance of many existing methods on intricate, real-world field shapes under strict time constraints is still understudied.
\subsection{Contributions}
\label{sec:contrib}

The primary novelty of this work lies in the synthesis and evaluation of an end-to-end operational framework for agricultural crop spraying. While individual components such as $p$-Median siting, $k$-means clustering, and \ac{TSP}formulations are established in the literature, there remains a gap in understanding how these upstream decisions cascade through a multi-stage planning pipeline to impact downstream performance. By transitioning these theoretical models into an applied benchmarking framework, we systematically quantify the trade-offs between computational efficiency and operational coverage. The specific contributions of this paper are split into four categories.

First, we introduce a novel application of the \ac{CoVRP} framework tailored specifically for precision crop spraying, departing from the rigid, traditional ``lawnmower'' path planning often found in agricultural optimization. This formulation provides a highly adaptable definition of coverage. While the baseline evaluation in this work utilizes a fixed-radius circular footprint for computational benchmarking, the underlying \ac{MIP} formulation is natively designed to accommodate complex, asymmetric coverage geometries allowing the optimization to inherently account for dynamic environmental factors without altering the core routing architecture.

Second, this paper delivers a controlled, large-scale comparison of facility siting strategies and their cascading impact on the full routing pipeline. We implement a load-balancing \ac{MIP} optimizer and a heuristic-driven, greedy $p$-Median method. By feeding the output of each into the same downstream routing algorithm, we provide a clean, head-to-head evaluation of how the initial facility placement choice affects later solution quality and planning latency. We demonstrate that using the $p$-Median heuristic reduces the runtime of facility siting by three orders of magnitude while only reducing total coverage by 4\%.

Third, we introduce and systematically evaluate a tiered decomposition method for solving large-scale \ac{CoVRP}s. The approach first partitions the operational area into smaller sub-regions, solves a high-level \ac{TSP} through the cluster centers to establish a macro-level tour, and finally computes detailed coverage routes within each sub-region before stitching them into a complete plan. By conducting an extensive parameter sweep, varying the number of clusters and distance between routing nodes, we map the performance-runtime frontiers. We identify that utilizing 6 to 8 clusters for each route leads to a significant reduction in runtime with minimal impact on coverage metrics.

Fourth, we establish the practical relevance and scalability of the integrated framework by testing it on a broad set of real and synthetic field geometries that vary in boundary complexity and target density. Recognizing the temporal constraints of real-world deployments, wall-clock runtime is reported as a primary metric alongside solution quality, demonstrating that this pipeline successfully achieves minute-scale planning on commodity hardware.
\section{Materials and Methods}

This section details the proposed operational framework for agricultural coverage routing. We first formalize the problem setting by defining the geometric environment, vehicle constraints, and the concept of spatial coverage. Following this mathematical foundation, we outline the end-to-end tiered solution architecture, detailing both the initial facility selection strategies and the subsequent multi-stage routing algorithm designed to generate feasible, coverage-maximizing flight paths under computational time limits.

\subsection{Problem Setting}
\label{sec:setting}

This section introduces the routing and siting problem in operational terms, along with mathematical notation to be used in the remainder of the paper.

\subsubsection{Field, Nodes, and Facilities}
A \emph{field} is the geographic region to be sprayed. To model this continuous area for optimization, we discretize it using two distinct sets of points: passive nodes and active nodes.

\emph{Active nodes}, denoted by the set $A$, represent a set of discrete, potential waypoints where a drone can execute a spraying operation. A drone's flight path, or \emph{route}, is constructed as a sequence of selected active nodes. Visiting an active node is the action that provides coverage to all passive nodes within its coverage area.

\emph{Passive nodes}, denoted by the set $P$, form a fine-grained grid that represents the total area requiring treatment. These nodes function as evaluation points to measure the quality of spray coverage. A passive node is considered ``covered" if it falls within the effective coverage area of a drone visiting a nearby active node. To aid in evaluation later, we also define 
\emph{Double Covered Passive Nodes}, denoted by the set $Q \subseteq P$, which are the passive nodes within the coverage neighborhood of multiple visited active nodes. The primary objective of the coverage routing algorithm is to maximize the number of unique passive nodes covered.

\begin{figure}
    \centering
    \includegraphics[width=0.5\linewidth]{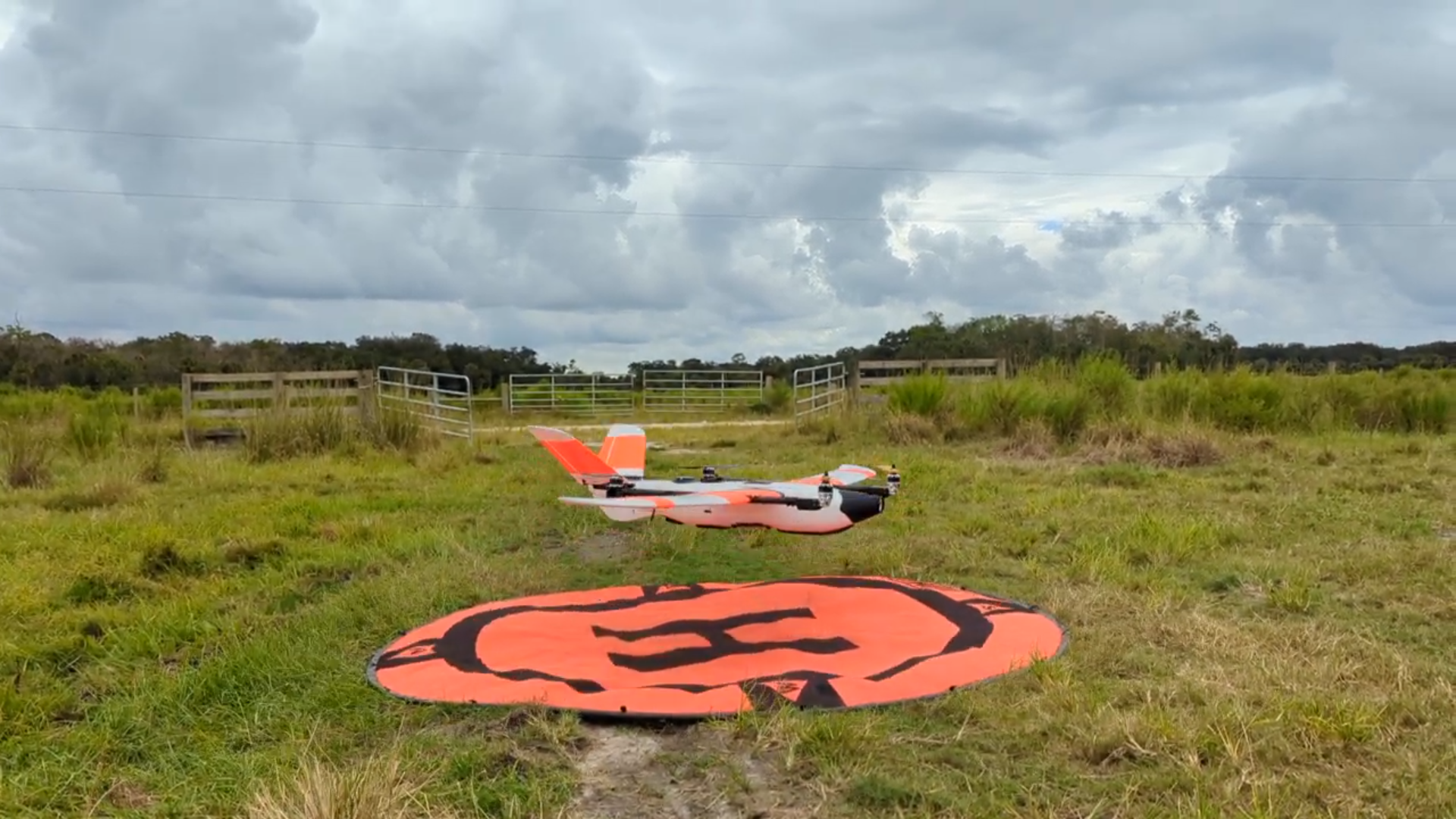}
    \caption{An example of a launch and recovery ``facility," which we approach from a facility siting optimization perspective.}
    \label{fig:dronetakeoff}
\end{figure}

Launch and recovery points for operators are called \emph{facilities}, an example of which is shown in Figure \ref{fig:dronetakeoff}. They are places where aircraft would take off, land, and swap batteries. We distinguish \emph{candidate facilities} $F^{\mathrm{cand}}$ from \emph{chosen facilities} $F\subseteq F^{\mathrm{cand}}$. Candidates are the acceptable locations surveyed before planning (e.g., field edges, roadsides, staging pads). The optimization selects a subset $F$ subject to a \emph{facility budget} $N_{open}$. Chosen facilities are the only depots from which routes start and end.

\subsubsection{Vehicles and Distances}
We plan for a homogeneous fleet of $K$ identical multirotor drones. The travel cost between any two locations $i$ and $j$ (which can be active nodes or facilities) is given by the distance $D_{ij}$. Each drone's operational range is limited by its battery capacity, which is simplified and modeled as a maximum flight distance, $D_{\text{max}}$. This parameter imposes a hard constraint on route construction: any single sortie, which must originate from and terminate at a facility, cannot have a total path length exceeding $D_{\text{max}}$.

The primary function of the vehicle is to provide coverage. When a drone visits an active node, it is assumed to cover any passive nodes that fall within its coverage footprint. In the nominal case, this footprint is modeled as a simple circle with a fixed spray radius, centered at the active node. This idealization provides a computationally tractable representation of spray dispersal. However, the framework is designed to accommodate more sophisticated models where the geometry of this footprint can be generalized. For example, the base circular area could be transformed by a function that accounts for external factors like wind, resulting in an elliptical or otherwise irregular coverage shape. Regardless of the specific model employed, a passive node is considered ``covered" if it is located within the boundary of the coverage footprint generated by a visited active node.

\subsubsection{Tiered solution overview}
Vehicle routing problems and their extensions are known to be NP-hard. As a consequence, the computation time to solve them grows rapidly as the problem size increases. To meet minute-scale planning, we solve the CoVRP in tiers: (1) A siting stage selects $F$; (2) A high-level routing stage orders clusters and assigns nodes to facilities; (3) Subproblems then solve coverage tours per cluster, then each subproblem is combined together into a complete route. This structure reduces instance size, enables warm starts, and exposes controllable parameters (cluster counts, arc sparsity, and time limits) whose effects we measure experimentally.

\begin{table}[htbp]
\centering
\caption{Nomenclature of sets, parameters, and decision variables.}
\label{tab:nomenclature}
\renewcommand{\arraystretch}{1.2} 
\begin{tabular}{l p{0.75\linewidth}}
\hline
\textbf{Symbol} & \textbf{Description} \\
\hline
\multicolumn{2}{l}{\textit{Sets and Indices}} \\
$A$ & Set of active nodes (potential drone waypoints) \\
$P$ & Set of passive nodes (evaluation points for field coverage) \\
$Q$ & Set of double-covered passive nodes ($Q \subseteq P$) \\
$F^{\mathrm{cand}}$ & Set of candidate launch and recovery facilities \\
$F$ & Set of chosen/activated facilities ($F \subseteq F^{\mathrm{cand}}$) \\
$I$ & Set of active routing nodes used in facility selection \\
$I'$ & Scaled-down, representative subset of active routing nodes ($I' \subset I$) \\
$J$ & Set of candidate facility locations in the optimization model \\
$S$ & Set of selected facilities during the $p$-Median heuristic \\
\hline
\multicolumn{2}{l}{\textit{Parameters}} \\
$K$ & Number of homogeneous multirotor drones in the fleet \\
$D_{ij}$ & Distance (e.g., Euclidean) between node $i$ and node $j$ \\
$D_{\text{max}}$ & Maximum flight distance per drone per battery cycle \\
$N_{\mathrm{open}}$ & Number of candidate facilities to select \\
$L_{\mathrm{ideal}}$ & Ideal balanced node load per facility ($|I| / N_{\mathrm{open}}$) \\
$\lambda$ & Weighting coefficient balancing distance and load/coverage objectives \\
$M$ & A large positive constant (Big-M) \\
$\alpha$ & Seed fraction for $p$-Median dispersion initialization \\
$r$ & Number of 1-swap refinement rounds in the $p$-Median heuristic \\
$w_m$ & Reward weight for covering passive node $m$ ($+1$ inside field, $-1$ outside) \\
$A_m$ & Pre-computed set of active vertices whose footprints cover passive node $m$ \\
$n$ & Total number of nodes in a given CoVRP subproblem tour \\
\hline
\multicolumn{2}{l}{\textit{Decision Variables}} \\
$x_j$ & Binary variable: $1$ if facility $j$ is selected, $0$ otherwise \\
$y_{ij}$ & Binary variable: $1$ if active node $i$ is assigned to facility $j$, $0$ otherwise \\
$z_j$ & Continuous variable: absolute load deviation for facility $j$ from $L_{\mathrm{ideal}}$ \\
$x_{ij}$ & Binary variable: $1$ if a drone travels directly from node $i$ to node $j$ \\
$y_i$ & Binary variable: $1$ if active node $i$ is visited on the tour, $0$ otherwise \\
$c_m$ & Binary variable: $1$ if passive node $m$ is covered, $0$ otherwise \\
$u_i$ & Auxiliary continuous variable used for subtour elimination \\
\hline
\end{tabular}
\end{table}
\subsection{Methodology}

\allowdisplaybreaks
\label{sec:methods}
The core of our methodology is a tiered optimization pipeline that first selects a set of launch and recovery sites (facilities) and then plans detailed coverage routes from each selected site. This section details the two alternative approaches for the initial facility selection stage: a baseline mixed-integer program designed for solution quality and a fast heuristic approach designed for scalability without sacrificing too much solution quality.

\subsubsection{Facility Selection}
Effective facility selection serves as a critical bridge between processing a raw operational field polygon and generating actionable coverage routes. By partitioning the total coverage area into smaller, facility-specific subproblems, the downstream routing becomes computationally tractable. To explore this stage, we implement two distinct siting methods: a load-balancing \ac{MIP} that co-optimizes facility workload and travel distance, and a fast $p$-Median heuristic that strictly minimizes assignment distance. While these two approaches optimize different mathematical objectives, the goal of this stage is not a direct algorithmic comparison. Instead, we contrast these methods to explicitly analyze how differing upstream siting objectives cascade through the tiered pipeline to affect downstream routing performance and final coverage quality (analyzed in Section 3.2).

\subsubsection*{Load-Balancing Optimization}
To establish a performance baseline, we first formulate the facility selection task as a \ac{MIP}. The primary objective of this model is to select a predetermined number of candidate facilities, $N_{\mathrm{open}}$, and assign the set of active routing nodes, $I$, to them. The objective function minimizes a sum of two operational goals: minimizing the total travel distance from active nodes to those facilities, and reducing the load deviation between the activated drone facilities. 

In this context, the ``load'' physically represents the total geographic area assigned to each launch and recovery site. By ensuring that this area is distributed evenly, the framework strategically limits the subsequent strain on drone battery capacity within any single sub-region. Because this initial stage of the tiered pipeline is focused on field partitioning and depot assignment, vehicle battery endurance is not included as a constraint here. Instead, endurance and flight distance bounds are strictly enforced during the route-generation phase detailed in Section 2.2.2.

To maintain computational tractability, the optimization is performed on a scaled-down, representative subset of the active nodes, $I' \subset I$. Once the optimal facility locations are determined using this subset, all active nodes in the full set $I$ are assigned to their nearest selected facility in a deterministic post-processing step. The formal optimization problem is defined as follows:

Let $J$ be the set of candidate facility locations. We define binary decision variables $x_j$ for each $j \in J$, where $x_j=1$ if facility $j$ is selected and 0 otherwise. Binary variables $y_{ij}$ for $i \in I', j \in J$ indicate if scaled-down active node $i$ is assigned to facility $j$. Finally, the continuous variables $z_j \ge 0$ mathematically represent the absolute deviation between the actual number of nodes assigned to facility $j$ and the ideal, perfectly balanced load, defined as $L_{\mathrm{ideal}} = |I| / N_{\mathrm{open}}$. The model is then:

\begin{align}
\min_{\,x,\,y,\,z} \quad 
  & \sum_{j} z_{j} \;+\; \lambda \sum_{i,j} D_{ij}\,y_{ij} \\[6pt]
\text{s.t.} \quad
  & \sum_{j} x_{j} = N_{\mathrm{open}},  \\[4pt]
  & \sum_{j} y_{ij} = 1, 
    \quad \forall\,i, \\[4pt]
  & y_{ij} \le x_{j}, 
    \quad \forall\,i,\,j, \\[4pt]
  & x_{k} \;+\; y_{ij} \;\le\; 1, 
    \quad \forall\,i,j,k:\; d_{k i} < d_{j i}, \\[4pt]
  & z_{j} \ge \sum_{i} y_{ij} \;-\; L_{\mathrm{ideal}}, 
    \quad \forall\,j,  \\[4pt]
  & z_{j} \ge L_{\mathrm{ideal}} \;-\; \sum_{i} y_{ij}, 
    \quad \forall\,j,  \\[4pt]
  & z_{j} \le M\,x_{j}, 
    \quad \forall\,j,  \\[4pt]
  & x_{j} \in \{0,1\}, 
    \quad y_{ij} \in \{0,1\}, 
    \quad z_{j} \ge 0. 
\label{eq:bal_flp}
\end{align}
where $D_{ij}$ is the Euclidean distance between node $i$ and facility $j$, and $\lambda$ is a weighting coefficient that balances the load deviation and distance costs. Constraint (2) ensures exactly $N^\mathrm{open}$ facilities are chosen. Constraint (3) enforces that each scaled-down node is assigned to exactly one facility, while (4) links assignments to open facilities. Constraint (5) defines the load deviation. While this MIP formulation produces well-balanced and spatially compact assignments, its computational expense grows rapidly with the number of candidate facilities and active nodes, making it a bottleneck for large-scale problem instances.

\subsubsection*{Hybrid $p$-Median Facility Selection}

The $p$-Median heuristic provides a computationally efficient method for facility selection by strictly minimizing the total assignment distance between active demand nodes and their nearest open facility. Given a set of candidate facilities $F$, active demand nodes $A$, and a distance matrix $D=\{D_{i,j}\}$ representing the distance from candidate facility $j\in F$ to demand node $i\in A$, the objective is defined as:
\[
\min_{\substack{S\subseteq F\\ |S|=N_{open}}} \quad \sum_{i\in A}\min_{j\in S} D_{i,j}.
\]

To solve this, the heuristic operates in three sequential stages, as detailed in Algorithm 2. The process begins with dispersion seeding, where a farthest-first initialization selects a fraction of the $p$ target facilities to ensure spatial distribution across the field and prevent initial clustering. Next, during the greedy completion stage, the algorithm iteratively evaluates unselected candidates and adds the facility that provides the largest reduction in total assignment distance; if no candidate yields an improvement, the remaining slots are filled via dispersion to avoid stalling. This is followed by an optional refinement phase, which utilizes a local search to evaluate 1-swap moves that exchange one selected facility for an unselected candidate, accepting the swap only if it strictly reduces the total assignment distance. Finally, upon completion of these stages, each demand node is deterministically assigned to its nearest open facility.

\begin{algorithm}
\caption{Greedy $p$-Median with Dispersion Seeding and Optional 1-Swap}
\begin{algorithmic}[1]
\State \textbf{Input:} $F$ candidates, $A$ demands, distances $D_{j,i}$;
\State \quad \quad \quad target $N_{open}$, seed fraction $\alpha$, swap rounds $r$
\State \textbf{Output:} Set of selected facilities $S$

\State $S \gets \text{FarthestFirst}(F,\lceil \alpha N_{open}\rceil)$ \Comment{Stage 1: Dispersion Seeding}
\State $m_i \gets \min_{j\in S} D_{j,i}$ for all $i\in A$ \Comment{Track min distance for each demand point $i$}

\While{$|S|<N_{open}$} \Comment{Stage 2: Greedy Completion}
  \State $j^\star \gets \arg\max_{j\in F\setminus S}\ \sum_{i\in A}\big(m_i-\min\{m_i,D_{j,i}\}\big)$ \Comment{Find facility with max total distance reduction}
  \State $\Delta \gets \sum_{i\in A}\big(m_i-\min\{m_i,D_{j^\star,i}\}\big)$
  \If{$\Delta \le 0$} \textbf{break} \EndIf
  \State $S \gets S\cup\{j^\star\}$;
  \State $m_i \gets \min\{m_i,D_{j^\star,i}\}$ for all $i$
\EndWhile

\If{$|S|<N_{open}$} \Comment{If greedy search stopped early...}
  \State $S \gets \text{FarthestFirst}(F,N_{open};\text{init}=S)$ \Comment{...fill remaining $p$ slots with dispersed points}
\EndIf

\For{$t=1$ \textbf{to} $r$} \Comment{Stage 3: Refinement Phase}
  \State $(a,b) \gets \arg\max_{a\in S,\ b\in F\setminus S}\ \sum_{i\in A}\Big(m_i-\min\{\min_{j\in S\setminus\{a\}}D_{j,i},D_{b,i}\}\Big)$ \Comment{Find best swap $(a \in S, b \notin S)$}
  \State $m_{i, \text{new}} \gets \min\{\min_{j\in S\setminus\{a\}}D_{j,i},D_{b,i}\}$ for all $i$
  \If{$\sum m_{i, \text{new}} < \sum m_i$}
    \State $S\gets (S\setminus\{a\})\cup\{b\}$; \quad $m_i \gets m_{i, \text{new}}$ for all $i$ \Comment{Perform the swap}
  \Else 
    \textbf{break} \Comment{Stop if no improving swap is found}
  \EndIf
\EndFor
\State \textbf{return} $S$
\end{algorithmic}
\end{algorithm}

By maintaining a vector of the current minimum distances $m_i=\min_{j\in S}D_{j,i}$ for each demand node, the algorithm avoids recomputing full assignment costs during each step. This allows the marginal gain of each candidate to be evaluated in $O(|A|)$ time. Consequently, the greedy completion stage requires $O\!\left(\left(N_{open}-|S|\right)\,|F|\,|A|\right)$ operations, and a single 1-swap pass requires $O\!\left(|S|\left(|F|-|S|\right)\,|A|\right)$ operations. Edge cases where $N_{open}\le 0$ or $N_{open}\ge |F|$ immediately return none or all facilities, respectively.

\begin{figure}[htbp]
    \centering
    \begin{subfigure}[b]{0.48\textwidth}
        \centering
        \includegraphics[width=\textwidth]{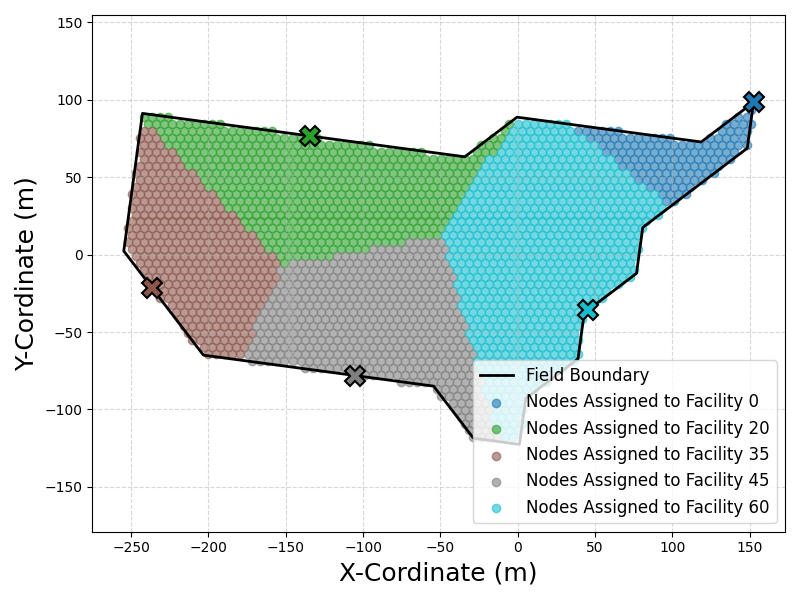}
        \caption{Imbalanced Assignment}
        \label{fig:facility-imbalanced}
    \end{subfigure}
    \hfill
    \begin{subfigure}[b]{0.48\textwidth}
        \centering
        \includegraphics[width=\textwidth]{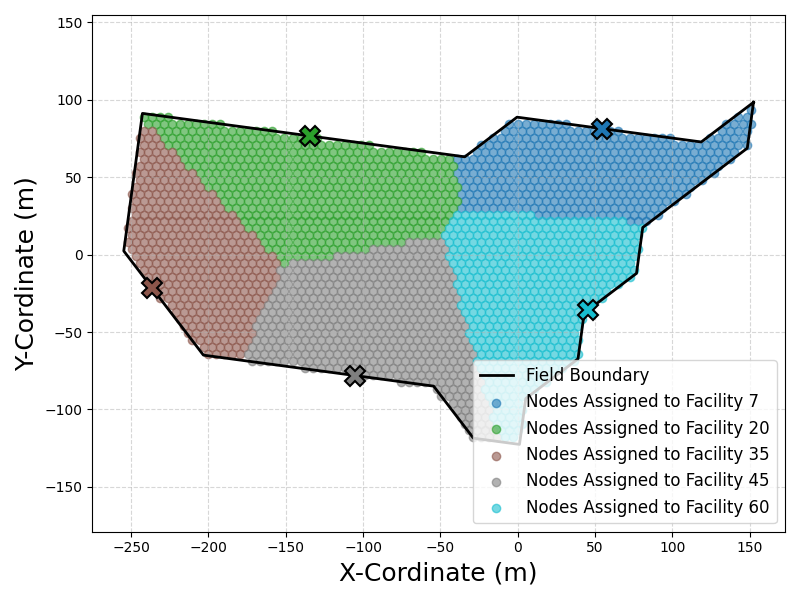}
        \caption{Balanced Assignment}
        \label{fig:facility-balanced}
    \end{subfigure}
    
    \caption{A geometric illustration of facility siting and active node assignments across an example field. (a) An unrefined assignment can result in imbalanced operational loads, where certain launch points (e.g., Facility 0 in dark blue) are assigned a disproportionate number of nodes. (b) A balanced assignment ensures a more equal distribution of nodes across all active facilities. This illustrates the motivation for incorporating balancing mechanisms, whether through MIP constraints or heuristic refinement swaps, to prevent vehicle endurance violations during the downstream routing phase.}
    \label{fig:facility-siting}
\end{figure}

\subsubsection{Tiered Coverage VRP}

\begin{figure}
    \centering
    \includegraphics[width=0.4\linewidth]{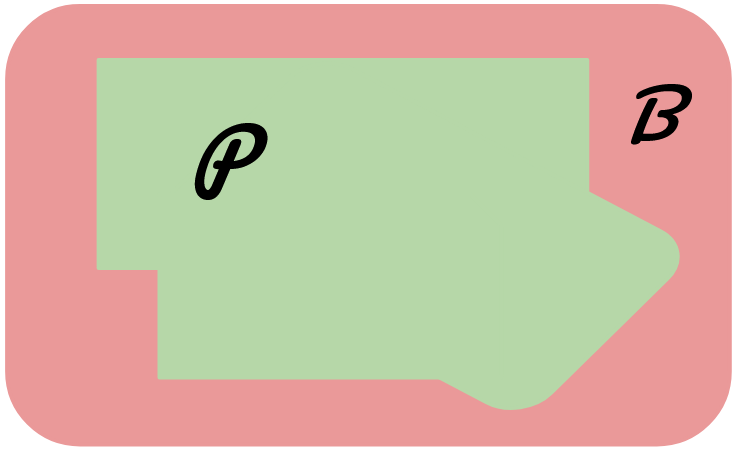}
    \caption{Example field illustrating how passive nodes are assigned coverage rewards. Nodes in the green region are assigned a reward of +1, while nodes in the red are assigned -1.}
    \label{fig:PassiveMockup}
\end{figure}

We begin with a field boundary polygon defined in GPS coordinates. To enable efficient spatial computation, these geographic coordinates are projected into a local planar coordinate system using the haversine function. The centroid of the field is translated to the origin $(0,0)$, and all other vertices are converted to meters relative to this point. Active routing nodes $A$ are generated via a hexagonal grid tiling clipped to the interior of the projected field polygon $\mathcal P$. Passive coverage nodes are laid on a finer square grid within a band surrounding the field's convex hull $\mathcal B$. Each passive node $p_m$ is assigned a reward

$$
  w_m = \begin{cases}
  +1, & p_m \in \mathcal P,\\
  -1, & p_m \in \mathcal B \setminus \mathcal P.
  \end{cases}
$$
This encourages solutions that remain inside the field and limit overspray outside of the field boundary (\cref{fig:PassiveMockup}). 

Once facilities are selected, route planning proceeds in two tiers. For each open facility, we perform k-means clustering on its assigned active nodes to identify high-level routing nodes, where high-level node 0 is the facility itself. The high-level tour is then computed by solving a standard Traveling Salesman Problem (TSP) utilizing the Miller-Tucker-Zemlin (MTZ) formulation for subtour elimination \cite{miller1960integer}. Each active node is reassigned to its nearest high-level node, forming groups. For each group, a CoVRP (Eq.~\eqref{eq:cvrp}) is solved independently, where endpoints are chosen based on TSP ordering to facilitate smooth concatenation. We formulate the CoVRP as the following mixed-integer program:

\begin{align}
\max_{x,y,c,u}\quad&
  \sum_{m=1}^{M} w_m\,c_m
  - \lambda\!\!\sum_{i\ne j}\!D_{ij}\,x_{ij}  \\
\text{s.t. }\;
&\sum_{i\ne j}\!D_{ij}\,x_{ij}\le D_{\max}, \\
&\sum_{j=1}^{n-1}x_{0j}=1,\;
 \sum_{i=0}^{n-2}x_{i,n-1}=1, \\
&\sum_{j\ne i}x_{ij}=y_i=
 \sum_{j\ne i}x_{ji},\quad(1\le i\le n-2),   \\  
&y_0=y_{n-1}=1, \\
&u_i-u_j+(n-1)x_{ij}\le n-2, &&\!(1\le i\ne j\le n-2), \\
&c_m\le\sum_{k\in A_m} y_k, &&\forall m, \\
&c_m\ge y_k\;(k\!\in\!A_m), &&\forall m, \\
&x_{ij}\!\in\!\{0,1\},\,y_i\!\in\!\{0,1\},\,
 c_m\!\in\!\{0,1\},\,u_i\!\in[1,n-1].
\label{eq:cvrp}
\end{align}
Here, binary variables $x_{ij}$ define the drone path, $y_i$ indicates whether a node is visited, $c_m$ denotes if passive node $m$ is covered, and $u_i$ supports subtour elimination. Each passive node $p_m$ has a pre-computed set $A_m$ of active vertices that can cover it. The objective balances coverage against distance, weighted by a user-defined penalty $\lambda$. Coverage can be defined in many different ways. In the simplest case, the spray can be modeled as a circular disk. In this case, any passive nodes within a defined radius of a visited active node are assumed to be covered. This can be adjusted to model other aspects, for instance, impacts on the spray due to wind.

\begin{figure}[htbp]
    \centering
    \begin{subfigure}[b]{0.48\linewidth}
        \centering
        \includegraphics[width=\linewidth]{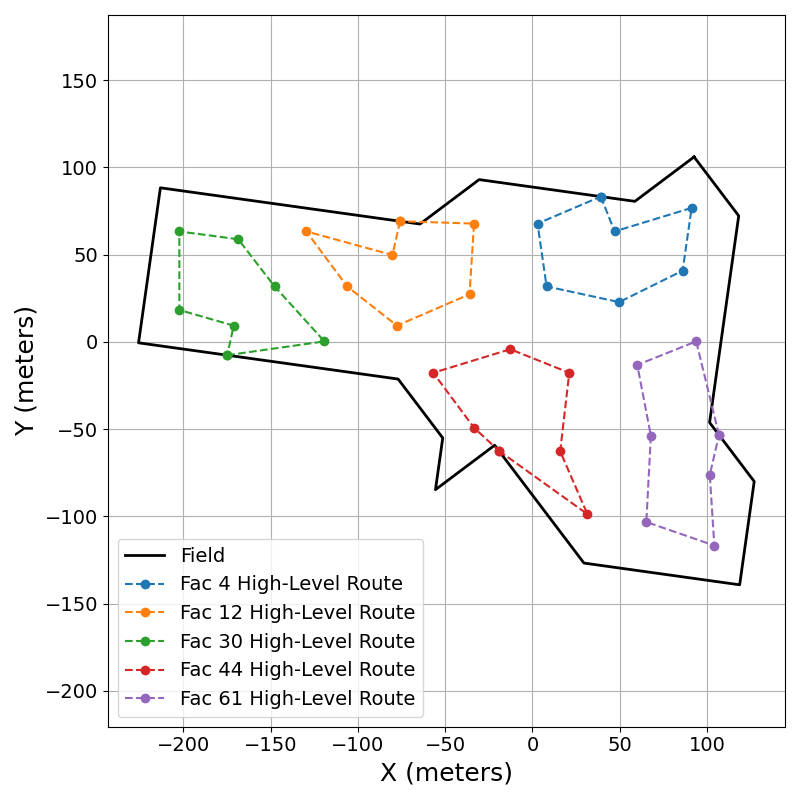}
        \caption{High-level TSP tour}
        \label{fig:route_tsp}
    \end{subfigure}
    \hfill
    \begin{subfigure}[b]{0.48\linewidth}
        \centering
        \includegraphics[width=\linewidth]{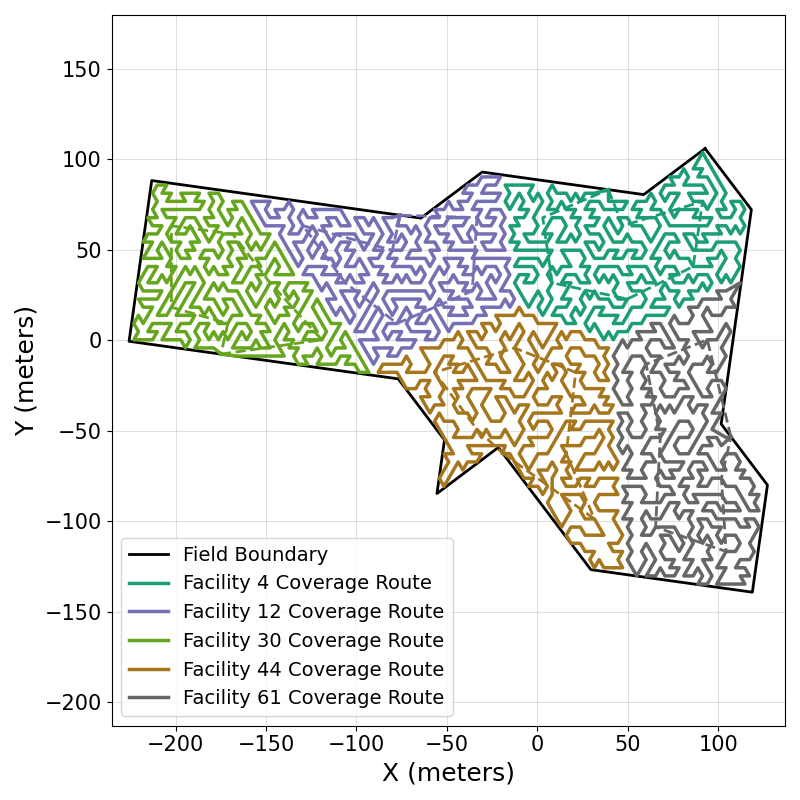}
        \caption{Complete coverage route}
        \label{fig:route_complete}
    \end{subfigure}
    \caption{An example of the tiered routing generation. (a) The macro-level tour established by solving a Traveling Salesman Problem through the cluster centers to sequence the sub-regions. (b) The final complete route after solving the \ac{CoVRP} for each cluster and applying 2-opt refinement, optimized to maximize coverage while adhering to distance constraints.}
    \label{fig:tiered_routing_example}
\end{figure}

After stitching the smaller sub-routes into a single overall route, this route is refined with the classical intra-route 2-opt heuristic \cite{croes1958twoopt,lin1965tsp}. The method selects two non-adjacent edges $(i,i+1)$ and $(j,j+1)$, removes them, reverses the sub-path between vertices $i+1$ and $j$, and reconnects the loose ends with the edges $(i,j)$ and $(i+1,j+1)$. Whenever the original pair of edges intersect, the replacement pair is strictly shorter by the triangle inequality, so repeated application systematically uncrosses every self-intersection introduced during subtour stitching. Because no vertices are added or removed, full coverage is preserved while the total route length is reduced without violating the global distance budget. The progression from the macro-level tour to the complete, fully refined route is illustrated in Figure \ref{fig:tiered_routing_example}.

\section{Results and Analysis}

This section presents the evaluation of the proposed agricultural coverage routing framework. To contextualize the performance outcomes, we first detail the experimental design, including the dataset properties, mission parameters, and the computational environment used for all simulations. Following this setup, we benchmark the facility selection methods—contrasting the optimization baseline against the $p$-Median heuristic—and subsequently analyze the impact of varying grid spacing and subgroup parameters within the tiered routing pipeline on both runtime and solution quality.

\subsection{Experimental Design}
\label{sec:experiment}
\subsubsection{Data set}

The agricultural field geometries used in this study were generated from binary raster maps representing predicted crop presence in Huron and Gratiot counties, Michigan, USA. The conversion of these raster images into the vector polygons required for route planning was accomplished using a computer vision pipeline detailed by Ni et al. \cite{ni2024machine}. In brief, the process utilizes contour detection to identify contiguous pixel clusters, followed by polygon simplification via the Douglas-Peucker algorithm to reduce vertex count while preserving the essential field shape. A final area thresholding step discarding polygons under 20 pixels is applied to remove noise and small artifacts. This pipeline yields a refined set of vector polygons that serve as the foundational data set for our experimental validation. For a comprehensive explanation of the algorithm's implementation and exact parameters, readers are directed to the work by Ni et al. 

The overall dataset contains around 1,500 individual fields, with many too large to be feasible for drone spraying. Contemporary agricultural multirotors are constrained by limited battery endurance and finite payload capacities. Consequently, covering large fields would require an impractical number of return-to-launch cycles and excessive battery swapping, falling outside the intended operational scope of continuous coverage routing. To ensure the robustness of the reported results across varying operational conditions, the experimental dataset comprised 20 distinct field instances. These fields ranged in size from 10 to 12 acres and were selected to encompass a diverse set of geometric boundaries, including long and narrow, roughly square, and non-convex polygons.

\subsubsection{Experimental Setup}

All simulations, models, and optimization algorithms were implemented in Python. The \texttt{gurobipy} library was utilized to interface the Python environment with the Gurobi solver for the \ac{MIP} formulations. All computational experiments were executed on a system equipped with an AMD Ryzen 5 5600X 6-core processor operating at a 4 GHz clock speed.

In all experiments, the nominal spray radius of the drone was fixed at 3 m. Each individual route was constrained by a distance budget of 2 km. This budget is derived from a total flight time of 15 minutes (900 seconds) at a constant speed of approximately 2.25 m/s ($\approx$ 5 mph), then rounding down to provide a small safety margin. 

To systematically estimate the required number of launch and recovery facilities for a given field, an upper bound for the total coverage area per drone sortie was calculated by multiplying the flight distance endurance by the drone's effective spray diameter. A 0.9 scaling factor was applied to this theoretical maximum to provide a margin of error for route inefficiencies and overlapping flight paths. The total field area was then divided by this scaled capacity to estimate the minimum number of discrete routes required for full coverage.

Passive nodes, representing the area to be covered, were generated on a square grid with a uniform spacing of 2 m. Active nodes, representing potential routing waypoints, were generated on an evenly spaced hexagonal grid; this structure was chosen for the more computationally intensive active nodes as it is the most efficient way to tile a plane. A set of 72 candidate facility locations was generated by spacing them evenly on the field boundary, positioned radially from the centroid of the field.

\subsubsection{Facility Selection Comparison}

The performance of the two distinct facility selection methods outlined in Section \ref{sec:methods} was evaluated. To investigate methods for reducing computational load, a scale-down strategy was tested for the load-balancing optimization, where the active node set was reduced by a factor of $n$ by including only every $n^{th}$ active node. For the $p$-Median heuristic, two core parameters were adjusted: the seed fraction (which dictates the initial number of facilities selected by dispersion) and the use of a final refinement swap round to improve the solution. During this phase of testing, all routing optimization parameters were held constant (5.7 m active node spacing, 6 subgroups per route) to isolate the impact of the facility selection method. 

\subsubsection{Routing Optimization Comparison}

Based on the results of the facility selection comparison, the $p$-Median heuristic was chosen for all subsequent routing tests. It was initialized with 2 seed facilities using the dispersion method, followed by 1 round of refinement swaps. To evaluate the tiered coverage vehicle routing algorithm, two key parameters were varied: the spacing in meters of the active nodes, and the number of subgroups used for each route.

\subsubsection{Evaluation Metrics}

The primary goal of this work is to decrease computational time to make the framework feasible for operational planning. Consequently, the run-time for both the facility selection and routing optimization phases was recorded as a key performance metric.

Solution quality was assessed using two coverage metrics. The passive coverage percentage, defined as the fraction of all passive nodes located within the coverage radius of a visited active node, served as a proxy for the total area sprayed. To measure redundant application, coverage efficiency was calculated as $1-|Q|/|P|$. This metric quantifies the degree of overlapping coverage, where a higher value indicates less ``wasted'' spray effort. The total distance of all combined routes was also tracked.

\subsection{Comparison Against Baselines}
To evaluate the potential advantages and limitations of the proposed framework, our results are benchmarked against two standard baseline methodologies prevalent in the coverage routing literature: an exact mathematical formulation and a non-decomposed routing approach.

First, following the comparative methodologies seen in facility location literature \cite{chauhan2019mcflpd}, we utilize the load-balancing \ac{MIP} formulation as the exact baseline. While MIPs guarantee optimal solutions, they represent the computational bottleneck typical of exact solvers. By benchmarking against this baseline, we explicitly quantify the trade-off of our proposed $p$-Median heuristic: a computational speedup of three orders of magnitude at the cost of a nominal 4\% degradation in passive coverage. 

Second, to validate the tiered routing architecture, we benchmark against the non-decomposed routing baseline, mirroring the evaluation strategy used by Cao et al. for complex environments \cite{cao2020hierarchical}. In our framework, the non-decomposed baseline is represented by the test cases where the number of subgroups is 1 (i.e., attempting to solve the \ac{CoVRP} over the entire assigned area). By comparing our tiered approach with an increasing number of subgroups against this single-group baseline, we demonstrate a runtime reduction from approximately 700 seconds to under 90 seconds, confirming an order-of-magnitude improvement over non-decomposed methods.

\subsection{Facility Siting: Runtime and Quality}

The performance of the optimized load-balancing and $p$-Median heuristic facility siting methods were evaluated, with key statistics presented in Tables \ref{tab:optimized_summary} and \ref{tab:pmedian_summary}.

\begin{figure}
    \centering
    \includegraphics[width=\linewidth]{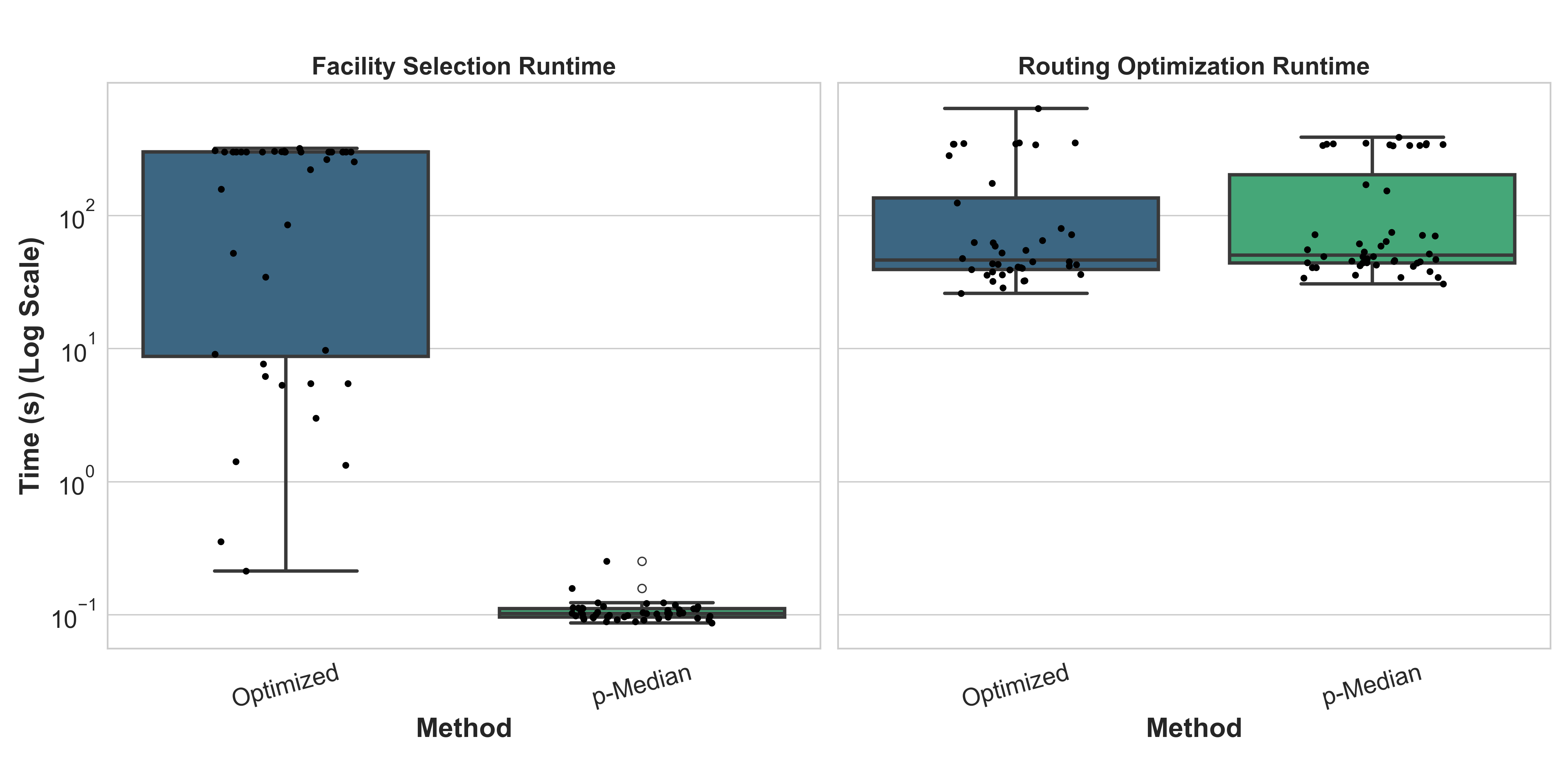}
    \caption{\textbf{Left}: Facility selection runtime vs. facility selection methods. \textbf{Right}: Routing optimization runtime vs. facility selection method.}
    \label{fig:fac_runtime_vs_siting_method}
\end{figure}

\begin{figure}
    \centering
    \includegraphics[width=0.7\linewidth]{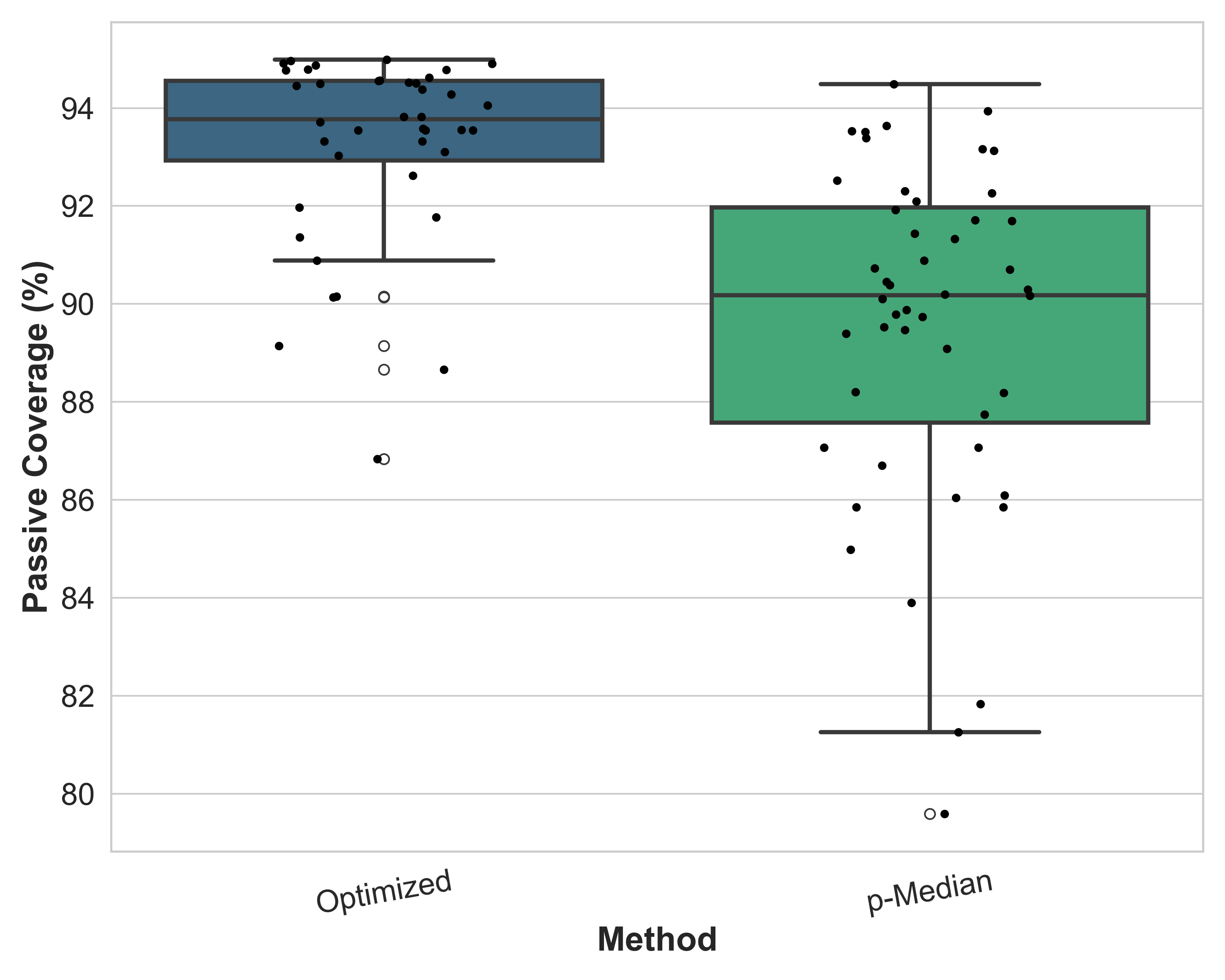}
    \caption{Passive node coverage percentage vs. facility selection method.}
    \label{fig:coverage_vs_siting_method}
\end{figure}

A significant difference in computational cost is immediately apparent. The p-Median heuristic is substantially faster, with facility selection times consistently under 1.2 seconds. In contrast, the optimized load-balancing method requires a much greater runtime, ranging from approximately 97 to 227 seconds (\cref{fig:fac_runtime_vs_siting_method}). This speed, however, comes with a trade-off in solution quality. The optimized method consistently achieves the highest passive coverage at around 93\%, while the p-Median heuristic shows a degradation in this metric, with coverage typically ranging from 86\% to 92\% (\cref{fig:coverage_vs_siting_method}). Interestingly, the routes generated from the facilities selected by the p-Median heuristic are, on average, shorter than those from the optimized method, which consistently exceed 9100 m.

\begin{table}[ht]
\centering
\caption{Optimized Method: Performance statistics vs. active node reduction. \textbf{Pas. Cov.} and \textbf{Cov. Eff.} are the passive coverage and coverage efficiency percentages respectively.}
\label{tab:optimized_summary}
\sisetup{separate-uncertainty, table-align-text-post=false}
\begin{tabular}{l S[table-format=3.2(3.2)] S[table-format=2.1(1.1)] S[table-format=2.1(1.1)] S[table-format=4.1(4.1)]}
\toprule
\multirow{2}{*}{\textbf{\shortstack{Node\\Reduc.}}} & {\textbf{Select Time (s)}} & {\textbf{Pas. Cov (\%)}} & {\textbf{Cov. Eff. (\%)}} & {\textbf{Route Len. (m)}} \\
& {Mean $\pm$ Std Dev} & {Mean $\pm$ Std Dev} & {Mean $\pm$ Std Dev} & {Mean $\pm$ Std Dev} \\
\midrule
1.0  & 209.28 \pm 127.67 & 92.8 \pm 2.2 & 96.1 \pm 0.5 & 9113.3 \pm 1103.6 \\
5.0  & 189.80 \pm 153.32 & 93.6 \pm 1.6 & 96.0 \pm 0.4 & 9203.0 \pm 1028.7 \\
10.0 & 211.10 \pm 130.26 & 93.0 \pm 2.1 & 96.0 \pm 0.5 & 9144.8 \pm 1161.6 \\
20.0 & 226.78 \pm 136.57 & 93.2 \pm 2.6 & 96.0 \pm 0.5 & 9178.2 \pm 1250.4 \\
40.0 & 97.38 \pm 136.24  & 93.5 \pm 1.6 & 95.9 \pm 0.4 & 9182.3 \pm 984.2 \\
\bottomrule
\end{tabular}
\end{table}

\begin{table}[ht]
\centering
\caption{$p$-Median Method: Performance Statistics by Seed Fraction and Refine Swaps. \textbf{Pas. Cov.} and \textbf{Cov. Eff.} are the passive coverage and coverage efficiency percentages respectively.}
\label{tab:pmedian_summary}
\sisetup{separate-uncertainty, table-align-text-post=false}
\begin{tabular}{cc S[table-format=1.2(1.2)] S[table-format=2.1(1.1)] S[table-format=2.1(1.1)] S[table-format=4.1(4.1)]}
\toprule
\textbf{Seed} & \textbf{Swaps} & {\textbf{Select Time (s)}} & {\textbf{Pas. Cov (\%)}} & {\textbf{Cov. Eff. (\%)}} & {\textbf{Route Len. (m)}} \\
& & {Mean $\pm$ Std Dev} & {Mean $\pm$ Std Dev} & {Mean $\pm$ Std Dev} & {Mean $\pm$ Std Dev} \\
\midrule
\multirow{2}{*}{0.2} & 0 & 0.10 \pm 0.01 & 89.4 \pm 1.9 & 96.7 \pm 0.3 & 8646.1 \pm 956.5 \\
                     & 1 & 0.11 \pm 0.01 & 92.6 \pm 1.7 & 96.2 \pm 0.3 & 9057.2 \pm 973.5 \\
\midrule
\multirow{2}{*}{0.4} & 0 & 0.10 \pm 0.01 & 86.5 \pm 4.1 & 97.1 \pm 0.5 & 8359.9 \pm 1049.0 \\
                     & 1 & 0.11 \pm 0.01 & 90.3 \pm 2.6 & 96.6 \pm 0.4 & 8792.3 \pm 1100.2 \\
\midrule
\multirow{2}{*}{0.6} & 0 & 0.10 \pm 0.02 & 88.0 \pm 3.7 & 97.0 \pm 0.6 & 8507.1 \pm 1041.5 \\
                     & 1 & 0.13 \pm 0.05 & 90.3 \pm 3.2 & 96.6 \pm 0.5 & 8760.3 \pm 989.7 \\
\bottomrule
\end{tabular}
\end{table}

\subsection{Tiered Routing: Parameter Study}

\begin{figure}[h!]
    \centering
    \includegraphics[width=0.98\linewidth]{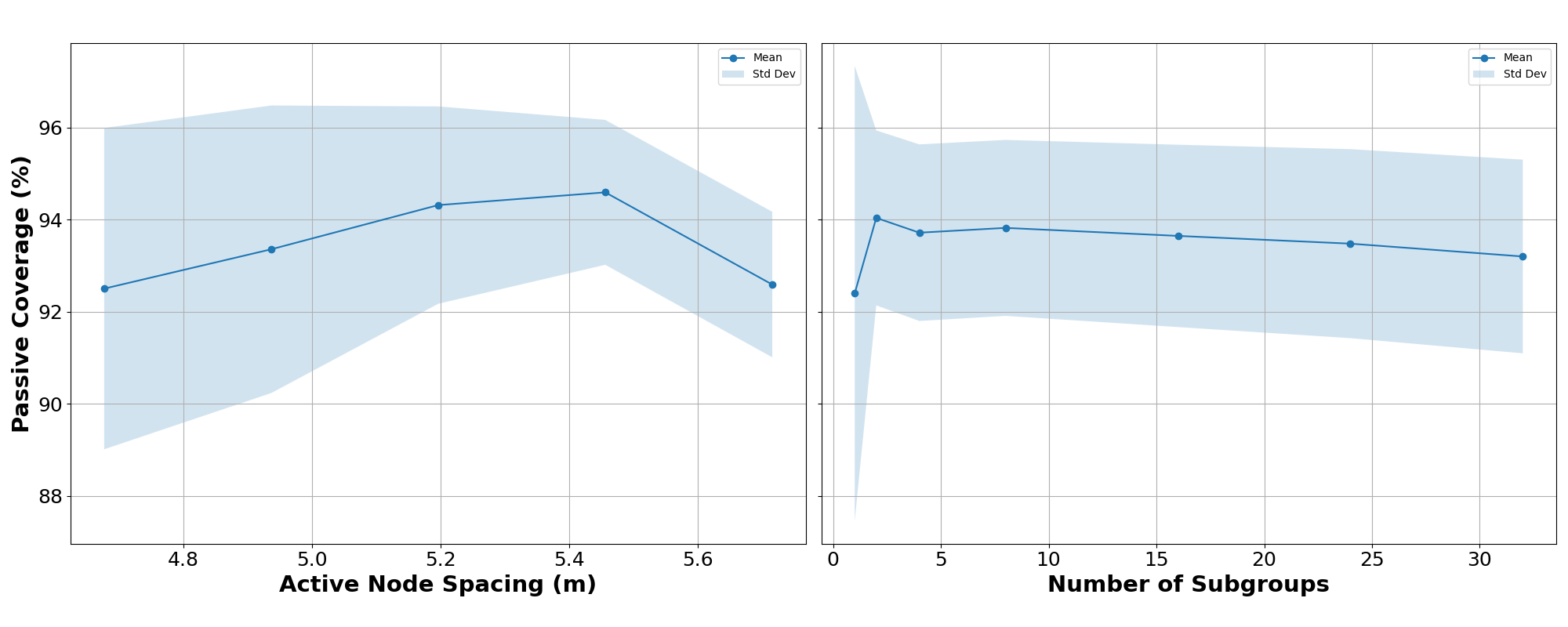}
    \caption{Passive coverage vs. active node spacing and number of subgroups.}
    \label{fig:passive_coverage_plot}
\end{figure}

\begin{figure}[h!]
    \centering
    \includegraphics[width=0.98\linewidth]{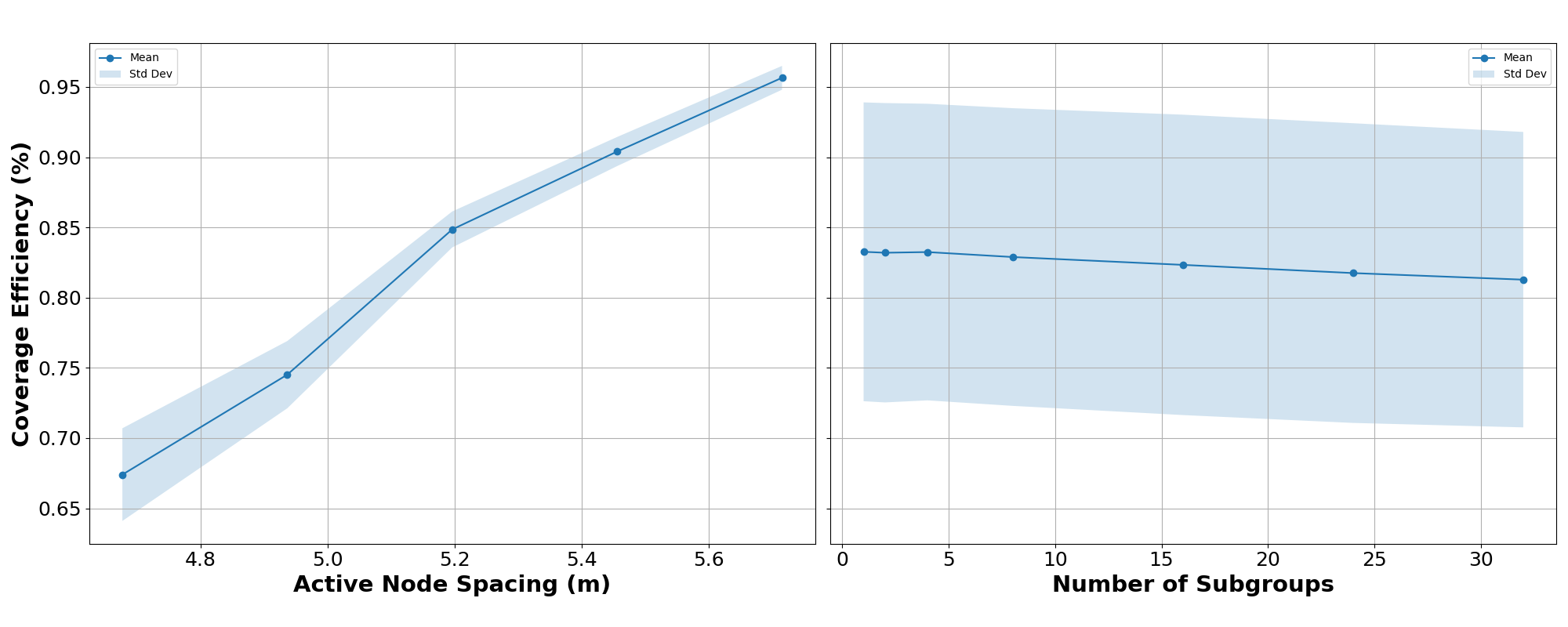}
    \caption{Coverage efficiency vs. active node spacing and number of subgroups.}
    \label{fig:coverage_efficiency_plot}
\end{figure}

\begin{figure}[h!]
    \centering
    \includegraphics[width=0.98\linewidth]{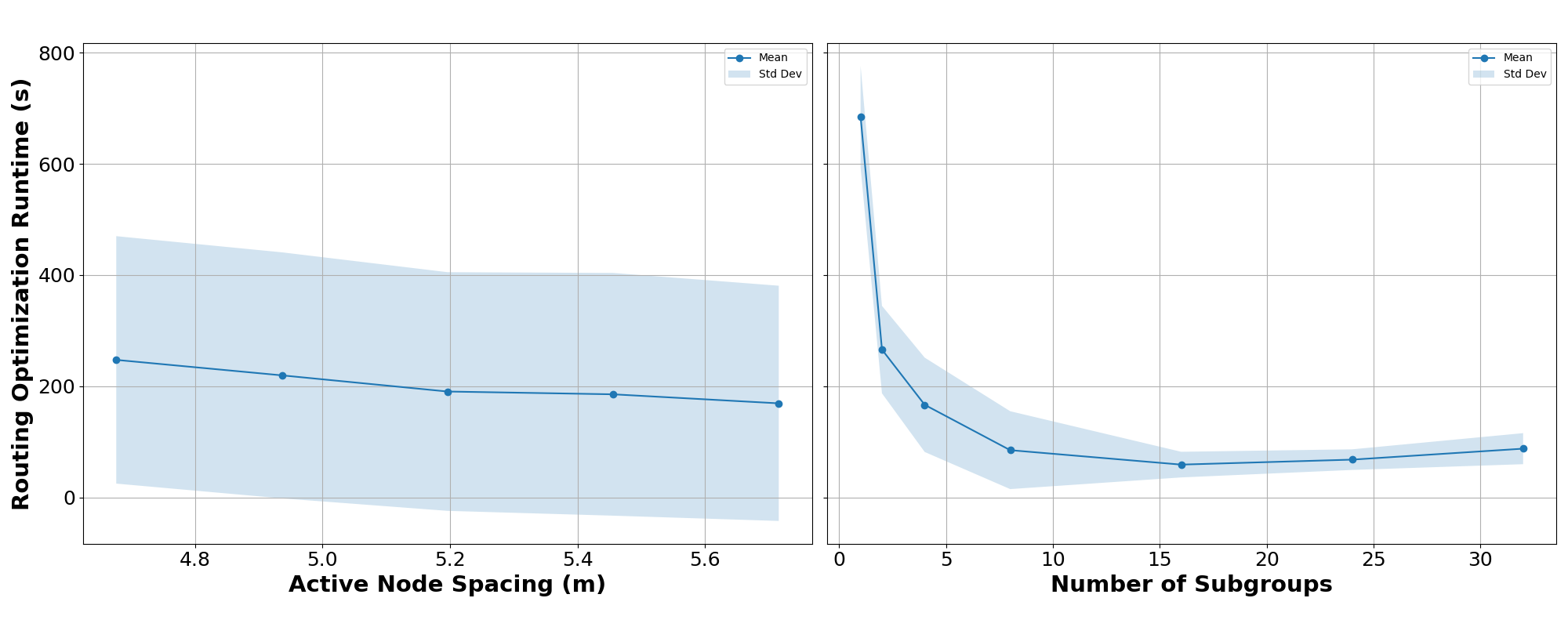}
    \caption{Routing runtime vs. active node spacing and number of subgroups.}
    \label{fig:runtime_plot}
\end{figure}

\begin{figure}[h!]
    \centering
    \includegraphics[width=0.7\linewidth]{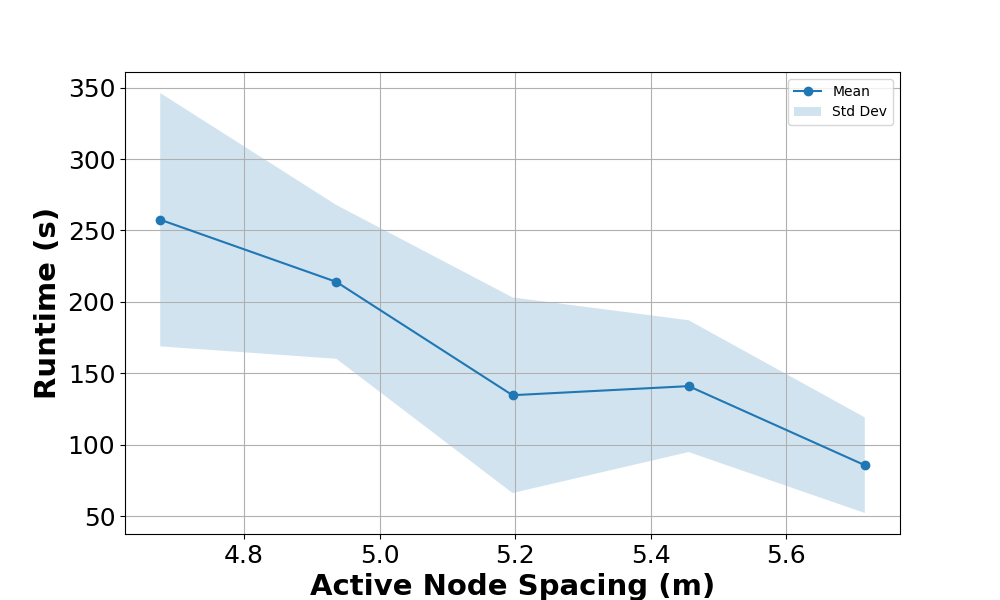}
    \caption{Routing runtime vs. active node spacing for a fixed number of 4 subgroups.}
    \label{fig:runtime_vs_spacing_plot}
\end{figure}

The tiered routing method's performance was analyzed by varying the effective grid spacing and the number of subgroups used for routing optimization. The results, summarized in Table \ref{tab:full_summary} and visualized in Figures \ref{fig:passive_coverage_plot}-\ref{fig:runtime_vs_spacing_plot}, reveal distinct trends.

As shown in Figure \ref{fig:runtime_plot}, the routing optimization runtime is strongly influenced by the number of subgroups. For all grid spacings, runtime decreases significantly as the number of subgroups increases from 1 to 16. This is expected, as dividing the problem into more sub-problems reduces the complexity of each optimization. However, beyond 16 subgroups, the runtime begins to level off or even slightly increase, suggesting that the overhead of managing a large number of sub-problems begins to outweigh the benefits of parallelization.

Coverage metrics also respond to parameter changes. Increasing the effective grid spacing generally leads to higher passive coverage and coverage efficiency, as seen in Figures \ref{fig:passive_coverage_plot} and \ref{fig:coverage_efficiency_plot}. This indicates that a coarser grid allows for more effective initial facility placement. Conversely, increasing the number of subgroups tends to cause a slight but consistent decrease in both passive coverage and efficiency, indicating a trade-off between runtime and solution quality.

\begin{table}[h!]
\centering
\caption{Performance Metrics (Mean $\pm$ Std. Dev.) vs. Grid Spacing and Number of Subgroups.}
\label{tab:full_summary}
\sisetup{
    separate-uncertainty,
    table-align-text-post=false
}
\begin{tabular}{
    c 
    S[table-format=3.2(2.2)] 
    S[table-format=2.2(1.2)] 
    S[table-format=1.2(1.2)] 
}
\toprule
\textbf{Subgroups} & {\textbf{Runtime (s)}} & {\textbf{Coverage (\%)}} & {\textbf{Efficiency (\%)}} \\
\midrule

\multicolumn{4}{c}{\textbf{Active Node Spacing: 4.68 m}} \\
\midrule
1  & 706.95 \pm 93.79 & 89.03 \pm 7.66 & 0.68 \pm 0.04 \\
2  & 351.02 \pm 56.70 & 93.71 \pm 1.69 & 0.68 \pm 0.03 \\
4  & 257.54 \pm 88.69 & 93.40 \pm 1.67 & 0.68 \pm 0.03 \\
8  & 161.90 \pm 119.89& 93.43 \pm 1.71 & 0.68 \pm 0.03 \\
16 & 83.85 \pm 28.91  & 92.97 \pm 1.83 & 0.67 \pm 0.03 \\
24 & 77.59 \pm 21.94  & 92.65 \pm 1.88 & 0.67 \pm 0.03 \\
32 & 93.06 \pm 31.18  & 92.32 \pm 1.97 & 0.66 \pm 0.03 \\
\midrule

\multicolumn{4}{c}{\textbf{Active Node Spacing: 4.94 m}} \\
\midrule
1  & 699.50 \pm 95.94 & 91.25 \pm 6.20 & 0.75 \pm 0.03 \\
2  & 311.78 \pm 66.71 & 94.25 \pm 2.18 & 0.75 \pm 0.03 \\
4  & 213.98 \pm 53.84 & 93.87 \pm 2.24 & 0.75 \pm 0.02 \\
8  & 90.89 \pm 35.94  & 93.95 \pm 2.21 & 0.75 \pm 0.02 \\
16 & 62.59 \pm 25.32  & 93.66 \pm 2.31 & 0.74 \pm 0.02 \\
24 & 70.10 \pm 20.59  & 93.43 \pm 2.38 & 0.74 \pm 0.02 \\
32 & 88.32 \pm 28.00  & 93.09 \pm 2.39 & 0.73 \pm 0.02 \\
\midrule

\multicolumn{4}{c}{\textbf{Active Node Spacing: 5.20 m}} \\
\midrule
1  & 678.22 \pm 90.02 & 94.07 \pm 3.03 & 0.85 \pm 0.01 \\
2  & 237.13 \pm 52.22 & 94.65 \pm 2.07 & 0.85 \pm 0.01 \\
4  & 134.66 \pm 68.48 & 94.33 \pm 2.06 & 0.86 \pm 0.01 \\
8  & 74.27 \pm 33.59  & 94.43 \pm 2.04 & 0.85 \pm 0.01 \\
16 & 53.04 \pm 12.78  & 94.41 \pm 2.08 & 0.85 \pm 0.01 \\
24 & 65.68 \pm 16.83  & 94.29 \pm 2.14 & 0.84 \pm 0.01 \\
32 & 90.47 \pm 31.88  & 94.04 \pm 2.25 & 0.84 \pm 0.01 \\
\midrule

\multicolumn{4}{c}{\textbf{Active Node Spacing: 5.46 m}} \\
\midrule
1  & 683.08 \pm 94.02 & 94.80 \pm 1.62 & 0.91 \pm 0.01 \\
2  & 229.01 \pm 59.35 & 94.76 \pm 1.62 & 0.91 \pm 0.01 \\
4  & 141.02 \pm 46.15 & 94.48 \pm 1.66 & 0.91 \pm 0.01 \\
8  & 48.84 \pm 11.57  & 94.67 \pm 1.62 & 0.91 \pm 0.01 \\
16 & 48.82 \pm 12.01  & 94.62 \pm 1.62 & 0.90 \pm 0.01 \\
24 & 62.95 \pm 15.75  & 94.52 \pm 1.72 & 0.90 \pm 0.01 \\
32 & 84.60 \pm 27.94  & 94.30 \pm 1.75 & 0.89 \pm 0.01 \\
\midrule

\multicolumn{4}{c}{\textbf{Active Node Spacing: 5.72 m}} \\
\midrule
1  & 657.39 \pm 99.60 & 92.84 \pm 1.64 & 0.96 \pm 0.01 \\
2  & 200.71 \pm 50.35 & 92.82 \pm 1.65 & 0.96 \pm 0.01 \\
4  & 85.67 \pm 33.47  & 92.52 \pm 1.67 & 0.96 \pm 0.01 \\
8  & 48.99 \pm 22.55  & 92.64 \pm 1.65 & 0.96 \pm 0.01 \\
16 & 46.63 \pm 11.23  & 92.57 \pm 1.63 & 0.95 \pm 0.01 \\
24 & 63.35 \pm 17.61  & 92.50 \pm 1.71 & 0.95 \pm 0.01 \\
32 & 81.88 \pm 26.48  & 92.25 \pm 1.69 & 0.94 \pm 0.01 \\
\bottomrule
\end{tabular}
\end{table}

\subsection{Discussion}
\label{sec:discussion}

This study provides guidance for practitioners by evaluating trade-offs between different facility siting methods, routing parameters, and environmental factors. When selecting a facility siting method, the $p$-Median heuristic appears to offer the most practical trade-off. It balances a significant reduction in facility selection runtime compared to the full optimization method while maintaining higher passive coverage than the Greedy MCLP approach. This suggests that for many applications where computational resources are a constraint, $p$-Median provides a robust compromise between speed and solution quality.

For tiered routing, the parameter study offers clear recommendations. The number of subgroups is a primary driver of runtime; a small number of subgroups provides a rapid decrease in computation time, but these gains diminish and eventually reverse as the number of subgroups becomes too large. Isolating for a fixed number of subgroups ($n = 4$), we observe that runtime trends downward as active node spacing increases. However, practitioners must be cautious, as increasing this spacing risks reducing passive coverage if nodes become too sparse for the spray radius to cover the intervening area. The relative independence of final solution quality (passive coverage and efficiency) from the number of subgroups is a key finding. It indicates that our method of decomposing the problem to accelerate computation does not inherently degrade the quality of the solution, which can instead be optimized through other parameters like grid spacing.

\subsection{Limitations of Work}
\label{sec:limitations}

While our proposed framework demonstrates strong performance, it is important to acknowledge its limitations, which present avenues for future work. The fidelity of the input data presents a challenge, as the data set was constructed from natural image data, and the automated processing introduced artifacts in some of the field polygons. These geometric inaccuracies can influence the coverage paths generated, potentially leading to suboptimal routing or incomplete coverage in real-world applications.

Methodologically, the overall optimization problem was decomposed into sequential stages: facility selection, subgroup partitioning, and the coverage VRP. Although this decoupling makes the problem tractable, these stages are inherently coupled, meaning decisions made in an early stage have cascading effects on the efficiency of the final routes. A monolithic or more integrated approach could theoretically yield globally optimal solutions, but at a much higher computational cost. Furthermore, the generated routes often contain sharp, non-kinematic turns that are not ideal for the dynamics of heavy quadrotor drones. This necessitates either reliance on the drone's on-board flight controller to approximate the path or a separate post-processing step to smooth the routes. While route smoothing could be incorporated directly into the optimization, this would introduce significant computational challenges and further constrain the path planning solution space.

Furthermore, the scope of our validation may limit generalizability. The algorithms were tested on fields ranging from 9 to 12 acres, which are relatively small compared to the 35-acre median field size in the full data set. The assumed drone parameters, such as a 5 mph speed and 15-minute flight time, push the limits of current commercial technology and necessitate numerous sorties even for these smaller fields, highlighting a potential bottleneck for larger-scale operations. Finally, the framework's performance is highly sensitive to several hyperparameters that require careful tuning, a process that can be complex and may require significant domain expertise when applying the model to new environments.
\section{Conclusions}
\label{sec:conclusion}

This study provides a comprehensive framework and practical guidance for optimizing agricultural drone operations by systematically evaluating the trade-offs inherent in different facility siting methods and routing parameters. Our analysis confirms that navigating these choices is critical for deploying efficient and effective autonomous systems in precision agriculture. The results show that the p-Median heuristic offers the most effective compromise for practitioners tasked with facility placement. It strikes a crucial balance, achieving a significant reduction in runtime compared to a full, computationally prohibitive optimization approach. This finding positions p-Median as a robust and practical tool for real-world planning scenarios.

Furthermore, the parameter study on tiered routing reveals that the number of subgroups is a primary and predictable driver of runtime. This allows for substantial computational savings through problem decomposition without inherently degrading the final solution quality. This key insight validates our sequential approach, demonstrating that operators can confidently accelerate planning by adjusting this parameter, while fine-tuning solution quality through other variables like grid spacing.

These findings also illuminate clear and promising avenues for future research that build directly upon the limitations of this work. A primary focus should be exploring more integrated optimization approaches that solve the facility location and vehicle routing problems in a coupled fashion. While our sequential method is computationally efficient, it risks settling in local optima because early-stage decisions are made without full knowledge of their downstream consequences. A monolithic or tightly integrated model, though more complex, could yield globally superior solutions. Lastly, validating this framework on a broader range of larger and more complex field geometries is crucial for confirming its scalability and robustness. Testing against non-convex fields, terrain with varying elevation, and environments with obstacles will be essential for transitioning this methodology from a theoretical model to a universally applicable operational tool.

\section*{Acknowledgments}

This project was supported by NASA under Award No. 80NSSC24PB252. We thank Jack R. Liddle, Benjamin L. Magocs, Jim Gardner, and Husni Idris for their support.

\section*{Hardware}

Tests were conducted on a custom personal computing system using a Ryzen 5 5600X 6-Core processor with a clock speed of 4 GHz.


\appendix

\printbibliography

\end{document}